# On Włodarczyk's Embedding Theorem

Jürgen Hausen

**Abstract**

We prove the following version of Włodarczyk's Embedding Theorem: Every normal complex algebraic $\mathbb{C}^*$-variety $Y$ admits an equivariant closed embedding into a toric prevariety $X$ on which $\mathbb{C}^*$ acts as a one-parameter-subgroup of the big torus $T \subset X$. If $Y$ is $\mathbb{Q}$-factorial, then $X$ may be chosen to be simplicial and of affine intersection.

## Introduction

One of the main results of [10] is that every normal ($\mathbb{Q}$-factorial, smooth) algebraic variety can be embedded into a (simplicial, smooth) toric prevariety. The purpose of this article is to prove the following equivariant and enhanced version of Włodarczyk's result (see Theorem 9.1):

**Theorem.** *Every normal complex algebraic $\mathbb{C}^*$-variety $Y$ admits an equivariant closed embedding into a toric prevariety $X$ on which $\mathbb{C}^*$ acts as a one-parameter-subgroup of the big torus $T \subset X$. If $Y$ is $\mathbb{Q}$-factorial, then $X$ may be chosen to be simplicial and of affine intersection.*

Here a toric prevariety $X$ with big torus $T$ is said to be of affine intersection if for any two maximal affine $T$-stable charts of $X$ their intersection is again affine (cf. [2]). The importance of this notion is demonstrated by the following application:

Every simplicial toric prevariety of affine intersection is of the form $U/\!/H$ where $U$ is an open toric subvariety of some $\mathbb{C}^n$ and $U \to U/\!/H$ is a geometric prequotient for the action of a subtorus $H \subset (\mathbb{C}^*)^n$ (cf. [1], Prop. 3.1 and [2], Section 8). So the above Theorem yields the following generalization of the construction of an equivariant affine cone over a projective $\mathbb{C}^*$-variety:

**Corollary.** *Every $\mathbb{Q}$-factorial algebraic $\mathbb{C}^*$-variety $Y$ is equivariantly isomorphic to a geometric quotient space $A/\!/H$, where $A$ is a quasi-affine variety with a regular action of a torus $H \times \mathbb{C}^*$.* □

As a further possible application, we indicate in Section 10, how an embedding $\imath$ of the above type of a complex variety $Y$ with an effective action of a torus $H$ into a toric prevariety $X$ of affine intersection can be used to obtain nice global resolutions of coherent $H$-sheaves $\mathcal{F}$ on $Y$ (see Theorem 10.1):

**Theorem.** *There exist $T$-stable Weil Divisors $D_1, \ldots, D_r$ on $X$ such that setting $\mathcal{E} := \bigoplus \imath^* \mathcal{O}_X(D_i)$ one obtains an exact sequence $\mathcal{E} \longrightarrow \mathcal{F} \longrightarrow 0$ of $H$-sheaves. If $X$ is smooth, then $\mathcal{E}$ is locally free.*

In particular, for the case of a locally factorial $Y$ and $H$ being trivial, the existence of an embedding of $Y$ into a smooth toric prevariety of affine intersection yields again the classical resolution result of Kleiman.

Having in mind the use of an embedding of a variety into a toric prevariety of affine intersection, one may ask if our first Theorem holds for any normal variety $Y$. The answer is negative; in [7] an explicit example of a normal surface is given that admits neither embeddings into toric prevarieties of affine intersection nor into simplicial ones.

The article is organized as follows: In Section 1 we present a reduction of simplicial toric prevarieties to such of affine intersection. In Sections 2 and 3 we give some existence statements on homogeneous functions, which are the crucial tool to perform Włodarczyk's constructions equivariantly.

Sections 4 and 5 are devoted to the equivariant reformulation of the notion of an embedding system of functions. In Sections 6 to 8 the existence of a family of embedding systems of functions is proved. Section 9 consists of the proof of the main result. Finally, in Section 10 we discuss the problem of resolving coherent sheaves.

As the reader will notice, we make in many cases use of the ideas provided in [10]. However, the equivariant situation and the reduction to affine intersection require sometimes a different and more involved treatment, hence, to be complete we give the full proofs in our setting.

The author wants to thank M. Brion and also J. Włodarczyk for various helpful discussions.

## 1  A Construction for Simplicial Toric Prevarieties

Let $N$ be a lattice and let $N_\mathbb{R} := \mathbb{R} \otimes_\mathbb{Z} N$ denote the associated real vector space. For any subset $B$ of $N_\mathbb{R}$, we define $\operatorname{cone}(B)$ to be the convex cone generated by $B$. By a *cone in $N$* we mean a convex polyhedral rational cone $\sigma \subset N_\mathbb{R}$.

For a cone $\sigma$ in $N$ we denote its relative interior by $\sigma^\circ$ and if $\tau$ is a face of $\sigma$, then we write $\tau \prec \sigma$. A cone $\sigma$ in $N$ is called *simplicial* if it has a linear independent set of generators. For a strictly convex cone $\sigma$, we denote by $\sigma^{(1)}$ the set of its extremal rays.

A *fan* in $N$ is a finite set $\Delta$ of strictly convex cones in $N$ such that $\sigma, \sigma' \in \Delta$ implies $\sigma \cap \sigma' \prec \sigma$ and $\sigma \in \Delta$ implies that also every face of $\sigma$ lies in $\Delta$. For a given cone $\sigma$ the set $\mathfrak{F}(\sigma)$ of its faces is a fan in $N$. A *system of fans* in $N$ is a finite family $\mathcal{S} := (\Delta_{ij})_{i,j \in I}$ of fans in $N$ such that

$$\Delta_{ij} = \Delta_{ji}, \qquad \Delta_{ij} \cap \Delta_{jk} \prec \Delta_{ik}$$

holds for any $i, j, k \in I$. Here $\Delta \prec \Delta'$ means that $\Delta$ is a subfan of $\Delta'$. A system $(\Delta_{ij})_{i,j \in I}$ of fans is called *affine*, if for every $i \in I$ the fan $\Delta_{ii}$ is the fan of faces of a single cone. A system $\mathcal{S}$ of fans in $N$ is called simplicial, if all its cones are simplicial.

Recall from [2], that every affine system $\mathcal{S}$ of fans defines a toric prevariety $X_\mathcal{S}$ and that every map of systems of fans defines a toric morphism. In fact the functor $\mathcal{S} \mapsto X_\mathcal{S}$ is an equivalence of categories. We call a toric prevariety simplicial, if it arises from a simplicial system of fans.

A toric prevariety $X$ with acting torus $T$ is said to be of *affine intersection*, if for any two maximal $T$-stable affine charts $U, U' \subset X$ the intersection $U \cap U'$ is again affine. If $X = X_\mathcal{S}$ with an affine system $\mathcal{S} = (\Delta_{ij})_{i,j \in I}$ of fans, then $X$ is of affine intersection if and only if every $\Delta_{ij}$ is the fan of faces of a single cone $\sigma_{ij}$.

In the sequel we indicate how to reduce a given simplicial toric prevariety to one of affine intersection. Given a simplicial affine system $\mathcal{S} = (\Delta_{ij})_{i,j \in I}$ of fans in $N$ we construct a new one as follows: for any two $i, j \in I$ let $|\Delta_{ij}| = \bigcup_{\sigma \in \Delta_{ij}} \sigma$ denote the support of $\Delta_{ij}$. Set

$$\widetilde{\sigma}_{ij} := \mathrm{conv}(|\Delta_{ij}|), \qquad \widetilde{\Delta}_{ij} := \mathfrak{F}(\widetilde{\sigma}_{ij}).$$

Then $\widetilde{\mathcal{S}} := (\widetilde{\Delta})_{i,j \in I}$ is an affine simplicial system of fans in $N$. Moreover the maps $\mathrm{id}_N$ and $\mathrm{id}_I$ define a map $(G, \mathfrak{g})$ of the systems $\mathcal{S}$ and $\widetilde{\mathcal{S}}$ of fans in the sense of [2] (see ibid., Lemma 5.6). For $i \in I$ and $\tau \in \Delta_{ii}$, let $X_{[\tau,i]}$ denote the associated affine chart of $X_\mathcal{S}$ (see [2], Section 2). As a consequence of the fibre formula [2], Proposition 3.5, we obtain:

**1.1 Proposition.** *$X_{\widetilde{\mathcal{S}}}$ is of affine intersection and the toric morphism $g \colon X_\mathcal{S} \to X_{\widetilde{\mathcal{S}}}$ associated to $(G, \mathfrak{g})$ has the following properties:*

  i) *for every $i \in I$ the restriction $g|_{X_{[\sigma_{ii},i]}}$ maps $X_{[\sigma_{ii},i]}$ isomorphically onto $X_{[\widetilde{\sigma}_{ii},i]}$.*

  ii) *for every $i \in I$ we have $g^{-1}(X_{[\widetilde{\sigma}_{ii},i]}) = \bigcup_{j \in I} X_{[\widetilde{\sigma}_{ij},j]}$,*

  iii) *$g$ is universal with respect to toric morphisms from $X_\mathcal{S}$ into toric prevarieties of affine intersection.*

  iv) *If $X_\mathcal{S}$ is smooth, then also $X_{\widetilde{\mathcal{S}}}$ is smooth.* □

## 2 Reminder on $\mathbb{C}^*$-Actions

In this section we recall some well-known basic facts on $\mathbb{C}^*$-varieties and also their proofs. Let $Y$ be an irreducible complex algebraic variety and assume that an algebraic torus $H$ acts on $Y$ by means of a regular map $H \times Y \to Y$.

The fixed point set of $H$ is denoted by $Y^H$, and, for a point $y \in Y$, we denote by $H_y$ its isotropy group. If $V \subset Y$ is a $H$-stable open set and $f$ is a regular function on $V$, then we call $f$ *homogeneous*, if there is a character $\chi \in \mathrm{X}(H)$ such that $f(h \cdot y) = \chi(h) f(y)$ holds for every $(h, y) \in H \times V$.

**2.1 Lemma.** *Assume that $Y$ is affine and let $y \in Y$. Then we have*

  i) *If $H \cdot y$ is non-trivial, then there is an $f \in \mathcal{O}(Y)$ homogeneous with respect to a non-trivial character such that $f(y) \neq 0$.*

  ii) *If $H \cdot y$ is closed, then for any $\chi \in \mathrm{X}(H)$ that is trivial along $H_y$, there is an $f \in \mathcal{O}(Y)$ homogeneous with respect to $\chi$ such that $f(y) \neq 0$.*

**Proof.** To verify i), take $y \neq y' \in H \cdot y$ and choose $\widetilde{f} \in \mathcal{O}(Y)$ with $\widetilde{f}(y) = 1$ and $\widetilde{f}(y') = 0$. Then, since $H$ acts regularly, $\widetilde{f} = \sum_{i=1}^r f_i$ with functions $f_i \in \mathcal{O}(Y)$ that are homogeneous with respect to pairwise different characters. Note that $f_i(y') = 0$ if and only if $f_i(y) = 0$. So there is an $i$ with $0 \neq f_i(y)$ and $f_i(y) \neq f_i(y')$. For this $i$, the function $f_i$ is as wanted.

We prove ii). The orbit map $t \mapsto t \cdot y$ induces an isomorphism $H/H_y \to H \cdot y$. In particular any given $\chi \in \mathrm{X}(H)$ that is trivial along $H_y$, defines a regular function $g_\chi$ on $H \cdot y$. Since $Y$ is affine, $g_\chi$ is the restriction of a function $f_\chi \in \mathcal{O}(Y)$.

Now, $f_\chi = \sum_{i=1}^r f_i$ with functions $f_i \in \mathcal{O}(Y)$ that are homogeneous with respect to pairwise different $\chi_i \in X(H)$. By uniqueness of such representations, we see that all $f_i$ homogeneous with respect to a character $\chi_i \neq \chi$ vanish along $H \cdot y$. Thus $g_\chi$ is also the restriction of an $f_i$ that is homogeneous with respect to $\chi$.  □

A rational function $f \in \mathbb{C}(Y)$ is called homogeneous, if its domain of definition $\mathrm{Def}(f)$ is $H$-stable and $f$ is homogeneous on $\mathrm{Def}(f)$ in the above sense. The following basic observation will be used to obtain rational homogeneous functions.

**2.2 Remark.** *If $V \subset Y$ is an affine $H$-stable open set and $f \in \mathcal{O}(V)$ is homogeneous, then $\mathrm{Def}(f)$ is $H$-stable and $f$ is a homogeneous rational function on $Y$.*  □

In the sequel, let $Y$ be affine. Moreover, let $p: Y \to Y /\!/ H$ denote the algebraic quotient, i.e., $p$ is the regular map defined by the inclusion $\mathcal{O}(Y)^H \subset \mathcal{O}(Y)$. Recall that every fibre $p^{-1}(z)$ contains a unique closed orbit and this orbit lies in the closure of any other orbit contained in $p^{-1}(z)$.

**2.3 Lemma.** *If the action of $H$ is non-trivial, then every fibre of $p$ contains non-trivial orbits.*

**Proof.** Assume that some $p^{-1}(p(y))$ contains only trivial orbits. Then $p^{-1}(p(y))$ consists in fact of a single point. By semicontinuity of the fibre dimension, this holds near $y$, which in turn implies $Y = Y^H$.  □

For the remainder of this section, let $H = \mathbb{C}^*$. Then the character group of $\mathbb{C}^*$ is isomorphic to $\mathbb{Z}$, and a function $f \in \mathcal{O}(Y)$ is homogeneous, if and only if it satisfies $f(t \cdot y) = t^k f(y)$ for some integer $k$ and every $(t, y) \in \mathbb{C}^* \times Y$. If so, then $k$ is called the *weight* of $f$. In this terminology, Lemma 2.1 yields

**2.4 Lemma.** *Let $Y_z := p^{-1}(z)$ be a fibre of the quotient map consisting of a single non-trivial $\mathbb{C}^*$-orbit. Then there exist homogeneous functions of positive weight and homogeneous functions of negative weight that do not vanish along $Y_z$.*  □

A *global source* (resp. *sink*) of the $\mathbb{C}^*$-action is a component $F$ of the fixed point set such that $\lim_{t \to 0} t \cdot y \in F$ (resp. $\lim_{t \to \infty} t \cdot y \in F$) holds for every $y$ in $Y$. If there is a global source (sink) then every homogenous function of negative (positive) weight is trivial.

A fibre $p^{-1}(z)$ is called *regular* if $z$ is a regular point of $Y /\!/ \mathbb{C}^*$, there are no singular points of $Y$ in $p^{-1}(z)$ and the differential of $p$ is surjective at any point of $p^{-1}(z)$. We recall the following general fact.

**2.5 Remark.** *There are precisely two possibilities for the quotient map $p: Y \to Y /\!/ \mathbb{C}^*$:*

  i) *Each regular fibre of $p$ consists of a single closed non-trivial $\mathbb{C}^*$-orbit.*

  ii) *There is a global source or a global sink of the $\mathbb{C}^*$-action.*

**Proof.** The statement being clear for trivial actions, we may assume that $\mathbb{C}^*$ acts non-trivially. Consider the image $A := p(Y^{\mathbb{C}^*})$. We show, that the two cases i) and ii) correspond to $A \neq Y /\!/ \mathbb{C}^*$ and $A = Y /\!/ \mathbb{C}^*$.

Let $A = Y /\!/ \mathbb{C}^*$. Then, for every $z \in Y /\!/ \mathbb{C}^*$, there is a unique fixed point $y(z) \in p^{-1}(z)$. Lemma 2.1 provides for every $y \in Y \setminus Y^{\mathbb{C}^*}$ a homogeneous function $f_y$ of non-trivial weight $k_y$ with $f(y) \ne 0$. We claim that any two of these weights are of the same sign. Otherwise, we found $y_1, y_2$ such that $f_{y_1}^{n_1} f_{y_2}^{n_2}$ is invariant for some positive integers $n_i$. Thus

$$f_{y_1}^{n_1}(y) f_{y_2}^{n_2}(y) = f_{y_1}^{n_1}(y(z)) f_{y_2}^{n_2}(y(z)) = 0$$

holds for any $y \in p^{-1}(z)$, $z \in Y /\!/ \mathbb{C}^*$. This contradicts irreducibility of $Y$. Thus all $k_y$ are of the same sign. If they are all positive, we obtain $\lim_{t \to 0} f_y(ty) = 0$ and hence $\lim_{t \to 0} ty = y(z)$ for any $y \in p^{-1}(z)$, $z \in Y /\!/ \mathbb{C}^*$. In other words, $Y^{\mathbb{C}^*}$ is a global source for $Y$. Analogously we see that if all $k_y$ are negative, then $Y^{\mathbb{C}^*}$ is a global sink for $Y$.

Now, let $A \ne Y /\!/ \mathbb{C}^*$. This means in particular, that $\dim(Y^{\mathbb{C}^*}) = \dim(A) < \dim(Y /\!/ \mathbb{C}^*)$. Moreover, the generic fibre of $p$ consists of a non-trivial closed orbit and thus $\dim(Y /\!/ \mathbb{C}^*) = \dim(Y) - 1$ holds. We have to show that for $y \in Y^{\mathbb{C}^*}$ the fibre $p^{-1}(p(y))$ is not regular.

So, let $y \in Y$ be a fixed point. By Lemmata 2.1 and 2.3, there is a homogeneous function $f_1 \in \mathcal{O}(Y)$ of non-trivial weight, that does not vanish on some $y_1 \in p^{-1}(p(y))$. For dimension reasons, $\mathbb{C}^*$ acts non-trivially on the zero set $N(Y; f_1) := \{x \in Y; \ f_1(x) = 0\}$.

Since $y \in N(Y; f_1)$, it follows from Lemma 2.3 that $p^{-1}(p(y)) \cap N(Y; f_1)$ contains a non-trivial $\mathbb{C}^*$-orbit different from $\mathbb{C}^* \cdot y_1$. Hence $p^{-1}(p(y))$ is either reducible or its dimension is strictly bigger than one. In both cases it is not a regular fibre. □

## 3   An Equivariant Moving Lemma

Let $Y$ be an affine normal (hence irreducible) complex algebraic $\mathbb{C}^*$-variety and let $p \colon Y \to Y /\!/ \mathbb{C}^*$ denote the associated quotient map. It is known that the generic fibre of $p$ is irreducible (see e.g. [3], Prop. 4). Thus we have

**3.1 Remark.** *There exists a non-empty open $\mathbb{C}^*$-stable $V \subset Y$ such that*

i) *the restriction $p|_V$ has only regular fibres,*

ii) *for any $z \in p(V)$, the fibre $Y_z := p^{-1}(z)$ is irreducible; in particular, $Y_z \cap V$ is dense in $Y_z$.* □

Now, let $V \subset Y \setminus Y^{\mathbb{C}^*}$ be an open set as in the above remark. The following existence result on homogeneous functions is crucial for the construction of a $\mathbb{C}^*$-equivariant e.s.f.:

**3.2 Lemma.** *Let $D \hookrightarrow Y \setminus V$ and $A \hookrightarrow V$ be closed $\mathbb{C}^*$-stable subsets, where $D$ is of pure codimension one in $Y$. If $y \in Y \setminus D$, then there exists a homogeneous function $f \in \mathbb{C}[Y]$ such that*

i) $f|_D = 0$ *and* $f(y) \ne 0$,

ii) *the zero set $N(V; f) := \{x \in V; \ f(x) = 0\}$ contains no irreducible component of $A$.*

**Proof.** First we consider the case that $y$ is a fixed point of the $\mathbb{C}^*$-action. Since the quotient map $p$ separates disjoint closed $\mathbb{C}^*$-stable sets, one has $p(y) \notin p(D)$. Consequently, since $Y /\!/ \mathbb{C}^*$ is irreducible, the dimension of the (closed) set $p(D)$ is strictly less than $\dim(Y /\!/ \mathbb{C}^*)$.

We claim that $p(D) \cap p(V) = \emptyset$ holds. Otherwise, there were a point $z \in p(V) \cap p(D)$. By the properties 3.1 i) and ii) of $V$, one obtains

$$\dim(p^{-1}(z)) = \dim(p^{-1}(z) \cap V) = \dim(Y) - \dim(Y /\!/ \mathbb{C}^*)$$

Let $q \colon D \to p(D)$ denote the restriction of $p$ to $D$. Since $\dim(p(D)) < \dim(Y /\!/ \mathbb{C}^*)$, we have

$$\begin{aligned}
\dim(q^{-1}(z)) &\geq \dim(D) - \dim(p(D)) \\
&\geq (\dim(Y) - 1) - (\dim(Y /\!/ \mathbb{C}^*) - 1) \\
&= \dim(p^{-1}(z)).
\end{aligned}$$

Since $q^{-1}(z) \subset p^{-1}(z)$ holds and $p^{-1}(z)$ is irreducible by the choice of $V$, this implies $q^{-1}(z) = p^{-1}(z)$. In particular, there exist points in $D \cap V$, which is a contradiction. Hence our claim is verified.

As consequence of our claim, we see that for each irreducible component $A_i$ of $A$, the image $p(A_i)$ is not contained in $p(D)$. Thus we find a regular function $g \in \mathbb{C}[Y /\!/ \mathbb{C}^*]$ such that $g(p(y)) \neq 0$, $g|_{p(D)} = 0$ and $g$ does not vanish identically on any $p(A_i)$. Then $f := g \circ p$ is as wanted.

Now, assume that $y$ is not a fixed point. Consider the case $y \notin A$. Let $y' \in \mathbb{C}^* \cdot y$ be any point different from $y$. Choose a function $g' \in \mathbb{C}[Y]$ with $g'|_D = 0$, $g'(y) = 1$ and $g'(y') = 2$. Writing $g'$ as a sum of homogeneous functions, we see that readily some homogeneous part $g$ of $g'$ satisifies $g|_D = 0$ and $g(y') \neq g(y) \neq 0$. Let $k$ denote the weight of $g$. Since $g(y') \neq g(y)$ holds, $k$ is different from zero.

Let $A_1, \ldots, A_r$ be the irreducible components of $A$, ordered in such a manner that $A_1, \ldots, A_s$ are precisely those who lie in $N(V; g)$. Note that $A \subset V$ implies that no $A_i$ is of dimension zero. For each $i \leq s$, choose a point

$$y_i \in A_i \setminus \left(\overline{\mathbb{C}^* \cdot y} \cup \bigcup_{j \neq i} A_j\right).$$

Again by prescribing values, we find homogeneous functions $g_i \in \mathbb{C}[Y]$ of non-trivial weights $k_i$ such that

$$g_i|_D = 0, \quad g_i|_{A_j} = 0 \text{ if } j \neq i, \quad g_i|_{\overline{\mathbb{C}^* \cdot y}} = 0, \quad g_i(y_i) \neq 0.$$

If the quotient map $p \colon Y \to Y /\!/ \mathbb{C}^*$ satisfies to 2.5 ii), then the weights $k_i$ of $g_i$ have the same sign as $k$ and, with $m := k \prod k_i$, the desired function is given by

$$f := g^{m/k} + \sum_{i=1}^{s} g_i^{m/k_i}.$$

If we are in the situation of 2.5 i), then we find homogeneous functions $h_i \in \mathbb{C}[Y]$ of weights $l_i$ of the same sign as $k$ such that $h_i(y_i) \neq 0$. We may assume that $|l_i| > |k_i|$. Setting $n := |k| \prod |l_i + k_i|$, the desired function $f$ is obtained as

$$f := g^{n/|k|} + \sum_{i=1}^{s} (h_i g_i)^{n/|l_i + k_i|}.$$

If $y \in A_{i_0}$ for some $i_0$, then the desired function $f$ is obtained similarly as above: omit the function $g$ and choose $y_{i_0} = y$. $\square$

In general, there is no hope to extend Lemma 3.2 to actions of tori of higher dimension. We illustrate this by means of a two-dimensional torus action on a three-dimensional variety. Consider the following lattice vectors in $\mathbb{Z}^3$

$$v_1 := e_1 + e_3, \quad v_2 := e_2 - e_3, \quad v_3 := e_1 + e_2 + 3e_3, \quad v_4 := e_1 + e_2 - 3e_3.$$

Let $Y$ be the affine toric variety defined by $\sigma := \mathrm{cone}(v_1, \ldots, v_4)$ and denote by $H$ the subtorus of the acting torus $T$ of $Y$ that corresponds to the sublattice $(\mathbb{Z}^2 \times \{0\}) \subset \mathbb{Z}^3$. Then $Y /\!/ H$ is a point (see e.g. [1], Example 3.1).

Let $D_i$ be the irreducible $T$-stable Weil divisor on $Y$ corresponding to the ray $\varrho_i := \mathbb{R}_{\geq 0} v_i$ of $\sigma$. Set

$$D := D_3 \cup D_4, \quad V := Y_{\varrho_1} \cup Y_{\varrho_2} \setminus D, \quad A := (D_1 \cup D_2) \cap V.$$

Here $Y_{\varrho_i}$ denotes the open toric subvariety of $Y$ defined by the ray $\varrho_i \subset \sigma$. Note that $H$ acts freely on $Y$.

**3.3 Remark.** *Every $H$-homogeneous function $f \in \mathcal{O}(Y)$ with $f|_D = 0$ vanishes identically on some component of $A$.*

**Proof.** Let $f \in \mathcal{O}(Y)$ be homogeneous with respect to the action of $H$ such that $f|_D = 0$. Then $f$ is of the form $f = \sum_i \alpha_i \chi^{u_i}$ with certain characters $\chi^{u_i}$ of $T$. Note that regularity of $f$ implies $\langle u_i, v_j \rangle \geq 0$ for all $i, j$.

Since $f$ is homogeneous with respect to the action of $H$, each $u_i$ is of the form $(u_1, u_2, u(i))$ with fixed integers $u_1$ and $u_2$. Moreover, the condition $f|_D = 0$ means $\langle u_i, v_j \rangle > 0$ for any $i$ and $j = 3, 4$. In other words, we have

$$-(u_1 + u_2) < 3u(i) < u_1 + u_2$$

for all $i$. Now, assume that $f$ does not vanish identically on some component of $A$. Then there are indices $i, j$ such that $\chi^{u_i}|_{D_1} \neq 0$ and $\chi^{u_j}|_{D_2} \neq 0$, i.e., we obtain

$$\langle u_i, v_1 \rangle = \langle u_j, v_2 \rangle = 0.$$

This implies $u_1 + u_2 = 0$, which contradicts the above strict inequality. $\square$

## 4 Equivariant Systems of Functions

We reformulate here some basic ideas of [10] in an equivariant manner. Note that we will differ in details and in notation from [10]. For the convenience of the reader, also the proofs are given in our setting.

Let $Y$ be a normal complex variety and assume that an algebraic torus $H$ acts regularly on $Y$. For a Weil divisor $D = \sum n_i D_i$ on $Y$ we denote by $|D| := \bigcup_i D_i$ its support and, given a point $y \in Y$, by $D_y$ the divisor obtained by omitting from $D$ all components $n_i D_i$ with $y \notin D_i$.

Let $M$ be a finitely generated (multiplicative) subgroup of $\mathbb{C}(Y)^*$ such that each function $f \in M$ is homogeneous with respect to some character $\chi_f \in \mathrm{X}(H)$. Note that $M$ is a free $\mathbb{Z}$-module.

The set $M$ is called a *weak equivariant system of functions (weak equivariant s.f.)* on $Y$ if there is a finite set $\mathfrak{D}$ of $H$-stable effective Weil divisors on $Y$ such that the following properties hold:

**SF0)** If $D, D' \in \mathfrak{D}$ are different from each other, then they have no common components.

**SF1)** For every $f \in M$ the divisor $(f)$ is of the form $(f) = \sum_{D \in \mathfrak{D}} n_{D,f} D$.

**SF2)** For every $D \in \mathfrak{D}$, $p \in |D|$ and $q \in Y \setminus |D|$, there is an $f \in M$ with $(f)_p \geq D_p$ and $(f)_q \leq 0$.

In the sequel we assume that $M \subset \mathbb{C}(Y)$ is a weak equivariant s.f. on $Y$ and that $\mathfrak{D}$ is a finite set of $H$-stable effective Weil divisors on $Y$ satisfying to SF0 – SF2.

Note that by SF0 the coefficients $n_{D,f}$ in the representation of $(f)$ in SF1 are uniquely determined. This fact is used to associate to each $D \in \mathfrak{D}$ a lattice vector $\mathrm{val}_D \in N := \mathrm{Hom}(M, \mathbb{Z})$ as follows

$$\mathrm{val}_D \colon M \to \mathbb{Z}, \qquad f \mapsto n_{D,f},$$

where $n_{D,f}$ is the coefficient of $D$ in the representation of $(f)$ in SF1. We note some properties of the lattice vectors $\mathrm{val}_D$. For any point $y \in Y$ let

$$\mathfrak{D}^y := \{D \in \mathfrak{D}; \ y \in |D|\}.$$

**4.1 Remark.** *For every $f \in M$ one has:*

i) *For each $y \in Y$ there is a representation $(f)_y = \sum_{D \in \mathfrak{D}^y} \mathrm{val}_D(f) D$.*

ii) *If $D = \sum \alpha_i \widetilde{D}_i \in \mathfrak{D}$ with $\widetilde{D}_i$ irreducible, then $\alpha_i \mathrm{val}_D(f) = \mathrm{ord}_{\widetilde{D}_i}(f)$.* □

Now, let $N_\mathbb{R} := \mathbb{R} \otimes_\mathbb{Z} N$ be the real vector space associated to $N$. Then we define for every $y \in Y$ a convex polyhedral rational cone $\sigma_y \subset N_\mathbb{R}$ by

$$\sigma_y := \mathrm{cone}(\mathrm{val}_D; \ D \in \mathfrak{D}^y).$$

The *dual cone* of such a $\sigma_y$ is by definition $\sigma_y^\vee := \{u \in M_\mathbb{R}; \ u|_{\sigma_y} \geq 0\}$. Since the $\sigma_y$ are rational polyhedral cones, Remark 4.1 i) yields the following

**4.2 Remark.** *For each $y \in Y$ we have $\sigma_y^\vee = \mathrm{cone}(f \in M; \ (f)_y \geq 0)$.* □

This Remark states in particular, that the cones $\sigma_y$ do not depend on the choice of the set $\mathfrak{D}$. We collect further properties:

**4.3 Lemma.** *Every cone $\sigma_y$, $y \in Y$, is strictly convex.*

**Proof.** (cf. [10], Lemma 3.8). Let $y \in Y$. If $\mathfrak{D}^y$ is empty, then there is nothing to prove. If not, suppose

$$\sum_{D \in \mathfrak{D}^y} \alpha_D \mathrm{val}_D = 0$$

with non-negative $\alpha_D$. SF2 yields for any $D \in \mathfrak{D}^y$ a function $f_D \in M$ with $(f_D)_y \geq D_y$. Applying the above expression to such an $f_D$ yields $\alpha_D = 0$. □

**4.4 Lemma.** *Let $p, q \in Y$. Then $\sigma_p \cap \sigma_q$ equals $\mathrm{cone}(\mathrm{val}_D; \ D \in \mathfrak{D}^p \cap \mathfrak{D}^q)$ and is a common face of both, $\sigma_p$ and $\sigma_q$.*

**Proof.** (cf. [10], Lemma 3.9). For any $D \in \mathfrak{D}^p \setminus \mathfrak{D}^q$ choose, according to SF2, a function $f_D \in M$ with $(f_D)_p \geq D_p$ and $(f_D)_q \leq 0$. Analogously, choose for every $D \in \mathfrak{D}^q \setminus \mathfrak{D}^p$ a function $g_D \in M$ with $(f_D)_q \geq D_q$ and $(g_D)_p \leq 0$. Set

$$h := \left( \prod_{D \in \mathfrak{D}^p \setminus \mathfrak{D}^q} f_D \right) \left( \prod_{D \in \mathfrak{D}^q \setminus \mathfrak{D}^p} g_D^{-1} \right).$$

Then we have $h \in M$ and

$$\begin{aligned} \mathrm{val}_D(h) &> 0 \quad \text{if} \quad D \in \mathfrak{D}^p \setminus \mathfrak{D}^q, \\ \mathrm{val}_D(h) &= 0 \quad \text{if} \quad D \in \mathfrak{D}^p \cap \mathfrak{D}^q, \\ \mathrm{val}_D(h) &< 0 \quad \text{if} \quad D \in \mathfrak{D}^q \setminus \mathfrak{D}^p. \end{aligned}$$

Thus $h$ defines a linear function on $N$ that separates $\sigma_p$ and $\sigma_q$ and vanishes on $\mathrm{cone}(\mathrm{val}_D; D \in \mathfrak{D}^p \cap \mathfrak{D}^q)$. That proves the assertion. $\square$

We want to figure out situations, where the cones $\sigma_y$ are simplicial. This leads to the following notions:

The set $M$ is called an *equivariant system of functions (equivariant s.f.)* on $Y$ if there is a finite set $\mathfrak{D}$ of $H$-stable effective Weil divisors on $Y$ satisfiying SF0, SF1, SF2 as above and moreover:

**SF3)** For every $y \in Y$ and every homogeneous $g \in \mathbb{C}(Y)^*$ satisfying $(g)_y = \sum_{D \in \mathfrak{D}} n_D D_y$ there is a $f \in M$ with $(f)_y = n(g)_y$, for some positive $n \in \mathbb{Z}$.[1]

**4.5 Remark.** *Let $D \subset Y$ be an $H$-stable irreducible closed subvariety of codimension one and let $y \in D$. If $Y$ is $\mathbb{Q}$-factorial, then there exist an affine $H$-stable open neighbourhood $U \subset Y$ of $y$ and a homogeneous function $f \in \mathcal{O}(U)$ with $N(U; f) = D \cap U$.*

**Proof.** $\mathbb{Q}$-factoriality of $Y$ means that we find a neighbourhood $\widetilde{U}$ of $y$ and a function $\widetilde{f} \in \mathcal{O}(\widetilde{U})$ with $(\widetilde{f}) = n(D \cap \widetilde{U})$ for some $n \in \mathbb{Z}_{>0}$. We choose $\widetilde{U}$ and $\widetilde{f}$ in such a manner, that $n$ is minimal.

By Sumihiro's Theorem, $y$ has an affine $H$-stable neighbourhood $U$. By shrinking $\widetilde{U}$ appropriately, we can achieve that $\widetilde{U} = U_g$ for some $g \in \mathcal{O}(U)$. Then $\widetilde{f} = h/g^l$ with a function $h \in \mathcal{O}(U)$. Write $h = \sum h_i$ with homogeneous functions $h_i \in \mathcal{O}(U)$.

Note that each $h_i$ vanishes identically on $D \cap U$. Hence, by minimality of $n$, we have $(h_i)_y \geq (h)_y$. In other words, the function $h_i/h$ is regular at $y$. So, near $y$ we obtain

$$1 = \sum \frac{h_i}{h}.$$

Thus $(h_i)_y = (h)_y$ for some $i$. Hence, removing components of $(h_i)$ that do not contain $y$ from $U$ and making $U$ affine again by further shrinking, we obtain with $f := h_i$ the desired function. $\square$

**4.6 Lemma.** *Assume that $Y$ is $\mathbb{Q}$-factorial and $M$ is an equivariant s.f. on $Y$. Then, for each $y \in Y$, the set $\{\mathrm{val}_D; D \in \mathfrak{D}^y\}$ is linearly independent, i.e., the cone $\sigma_y$ is simplicial.*

---

[1] We give a weaker condition as Włodarczyk. The reason is that our construction will not guarantee his property in an equivariant manner.

**Proof.** (cf. [10], Lemma 3.11). Let $y \in Y$. Since $Y$ is $\mathbb{Q}$-factorial, Remark 4.5 provides for each $D \in \mathfrak{D}^y$ a homogeneous function $f_D \in \mathbb{C}(Y)^*$ with $(f_D)_y = nD_y$ for some positive $n \in \mathbb{Z}$. By SF3, we may assume that each $f_D$ belongs to $M$. Now, suppose that

$$\sum_{D \in \mathfrak{D}^y} \alpha_D \mathrm{val}_D = 0$$

holds for some real numbers $\alpha_D$, $D \in \mathfrak{D}^y$. Then we obtain for each fixed $D' \in \mathfrak{D}^y$ the identity

$$n\alpha_{D'} = \sum_{D \in \mathfrak{D}^y} \alpha_D \mathrm{val}_D(f_{D'}) = 0. \quad \square$$

Summing up the above considerations we obtain

**4.7 Proposition.** *If $M$ is an equivariant s.f. on $Y$, then $\Delta(M, Y) := \{\sigma_y;\ y \in Y\}$, is a fan in $N$. If $Y$ is $\mathbb{Q}$-factorial, then $\Delta(M, Y)$ is simplicial.* $\quad \square$

## 5 The Regular Map Arising from an S.F.

Let $H$ be an algebraic torus and let $Y$ be a normal complex algebraic $H$-variety. Moreover, let $M \subset \mathbb{C}(Y)$ be an equivariant s.f. with a set $\mathfrak{D}$ of $H$-stable effective Weil divisors on $Y$ satisfying to SF0 – SF3. For each $p \in Y$ we define

$$S_p := \sigma_p^\vee \cap M = \{f \in M;\ (f)_p \geq 0\},$$

$$U_p := \{y \in Y;\ \forall f \in S_p\ (f)_y \geq 0\} = \bigcap_{f \in S_p} \mathrm{Def}(f) \subset Y.$$

By Gordan's Lemma, each $S_p$ is a finitely generated semigroup and consequently each $U_p$ is an open subset of $Y$.

**5.1 Lemma.** *Let $p \in Y$ and $y \in U_p$. Then one has $\mathfrak{D}^y \subset \mathfrak{D}^p$, $\sigma_y \prec \sigma_p$, $S_p \subset S_y$ and $U_y \subset U_p$.*

**Proof.** Assume that there were $y \in U_p$ such that $\mathfrak{D}^y \not\subset \mathfrak{D}^p$. Then there is a $D \in \mathfrak{D}^y \setminus \mathfrak{D}^p$. By SF2), we find a function $f \in M$ with $f|_D = 0$ and $(f)_p \leq 0$. This means $(f^{-1})_p \geq 0$ and $y \notin \mathrm{Def}(f)$, a contradiction to the definition of $U_p$. Hence we obtain $\mathfrak{D}^y \subset \mathfrak{D}^p$.

So Lemma 4.4 yields that $\sigma_y$ is in fact a face of $\sigma_p$. Now, $S_p \subset S_y$ follows from the inclusion $\sigma_p^\vee \subset \sigma_y^\vee$ and $U_y \subset U_p$ follows by definition from $S_p \subset S_y$. $\quad \square$

Let $\mathbb{C}[S_p]$ denote the semigroup algebra determined by $S_p$. Moreover, let $X_{\sigma_p} := \mathrm{Spec}_\mathfrak{m}(\mathbb{C}[S_p])$ denote the affine toric variety associated to $S_p$ (see e.g. [6] ).

**5.2 Remark.** *For every $p \in Y$, the inclusion $S_p \subset \mathcal{O}(U_p)$ defines a homomorphism $\mathbb{C}[S_p] \to \mathcal{O}(U_p)$ and hence a regular map $\varphi_{M,p}\colon U_p \to X_{\sigma_p}$.* $\quad \square$

The homomorphism $M \to \mathrm{X}(H)$, $f \mapsto \chi_f$ defines an homomorphism $\psi$ from $H$ to the big torus $T := \mathrm{Hom}(M, \mathbb{C}^*)$ of $X_{\sigma_p}$. In particular $H$ acts on the toric variety $X_{\sigma_p}$ by means of $h*x := \psi(h) \cdot x$.

**5.3 Remark.** *Each map $\varphi_{M,p}\colon U_p \to X_{\sigma_p}$ is $H$-equivariant.*

**Proof.** For $f \in M$ let $\chi^f$ denote the associated character of $T$, i.e., the character with $\chi^f(t) = t(f)$. Then $\mathcal{O}(X_{\sigma_p})$ is generated as an algebra by the functions $\chi^f$, $f \in \sigma_p^\vee \cap M$. Moreover, for these $\chi^f$ we have

$$\chi^f(h*x) = \chi^f(\psi(h) \cdot x) = \chi_f(h)\chi^f(x).$$

In other words, each of the above $\chi^f$ is homogeneous with respect to the character $\chi_f \in \mathrm{X}(H)$. Since also every $f = \varphi_{M,p}^*(\chi^f)$ is homogeneous with respect to $\chi_f$ and $X_{\sigma_p}$ is affine, the claim follows. □

For $p \in Y$, let $\mathbb{C}\{S_p\}$ denote the subalgebra of $\mathcal{O}(U_p)$ that is generated by $S_p$. An equivariant system $M$ of functions on $Y$ is called an *equivariant embedding system of functions (equivariant e.s.f.)* on $Y$ if it has the following property

**SF4)** For every $p \in Y$ the set $U_p$ is affine and $\mathcal{O}(U_p) = \mathbb{C}\{S_p\}$.

**5.4 Lemma.** *Let $M \subset \mathbb{C}(Y)$ be an equivariant s.f. on $Y$ and let $p \in Y$. Then we have*

i) $\varphi_{M,y} = \varphi_{M,p}|_{U_y}$ *for every $y \in U_p$,*

ii) *If $M$ is an equivariant e.s.f., then $\varphi_{M,p}$ is a closed embedding.*

**Proof.** Assertion ii) is clear by definition. Assertion i) follows from the fact that for $y \in U_p$ the inclusion $S_p \subset S_y$ defines a commutative diagram

$$\begin{array}{ccc} \mathbb{C}[S_p] & \subset & \mathbb{C}[S_y] \\ \downarrow & & \downarrow \\ \mathcal{O}(U_p) & \subset & \mathcal{O}(U_y) \end{array}. \qquad \Box$$

Note that $H$ acts also on the toric variety $X_{\Delta(M,Y)}$ by means of $h*x := \varphi(h) \cdot x$. We conclude this section with stating the following equivariant version of [10], Proposition 3.5:

**5.5 Proposition.** *Let $M \subset \mathbb{C}(Y)$ be an equivariant s.f. on $Y$. Then the maps $\varphi_{M,p}$, $p \in Y$ glue together to an $H$-equivariant map $\varphi_{M,Y}\colon Y \to X_{\Delta(M,Y)}$. Moreover, one has*

i) $\varphi_{M,Y}^{-1}(X_{\sigma_p}) = U_p$ *for every $p \in Y$.*

ii) *If $M$ is an equivariant e.s.f., then $\varphi_{M,Y}$ is a closed embedding.*

**Proof.** The fact that the $\varphi_{M,p}$ glue together follows from Lemma 5.4 i). In order to verify i), we only have to show that $U_p$ contains $\varphi_{M,Y}^{-1}(X_{\sigma_p})$. This is done as follows:

$$\begin{aligned} \varphi_{M,Y}^{-1}(X_{\sigma_p}) &= \varphi_{M,Y}^{-1}\left(\bigcap_{f \in S_p} \mathrm{Def}(\chi^f)\right) \\ &= \bigcap_{f \in S_p} \varphi_{M,Y}^{-1}(\mathrm{Def}(\chi^f)) \\ &\subset \bigcap_{f \in S_p} \mathrm{Def}(\chi^f \circ \varphi_{M,Y}) \\ &= \bigcap_{f \in S_p} \mathrm{Def}(f) = U_p. \end{aligned}$$

Here $\chi^f$ denotes the character of $T$ corresponding to $f \in M$. To finish the proof, note that assertion ii) is a consequence of Lemma 5.4 ii) and Property i) which we just proved. □

# 6  Construction of a Family of $\mathbb{C}^*$-Equivariant E.S.F.

Let $Y$ be an arbitrary normal algebraic $\mathbb{C}^*$-variety. In view of Sumihiro's Theorem (see [9]), we find a finite cover of $Y$ by open $\mathbb{C}^*$-stable affine subvarieties $Y_i$, $i \in I$. For $i, j \in I$, let $Y_{ij} := Y_i \cap Y_j$.

In this section we construct a subgroup $M \subset \mathbb{C}(Y)$ defining an equivariant e.s.f. on each $Y_{ij}$ in four steps. The proofs of the non-trivial Steps (6.3 and 6.4) are given in the following section.

**6.1 Step.** *Choose a finite set $M^1 \subset \mathbb{C}(Y)$ of homogeneous functions such that for any two $i, j \in I$, the set $M^1$ contains generators of $\mathbb{C}[Y_{ij}]$.* □

Let $V \subset \bigcap_{i \in I} Y_i \setminus Y^{\mathbb{C}^*}$ be a non-empty open $\mathbb{C}^*$-stable affine subvariety such that every $f \in M^1$ is regular and non-vanishing on $V$, and moreover $V$ fullfills the conditions of 3.1 for all quotient maps $p_{ij} : Y_{ij} \to Y_{ij}/\!/\mathbb{C}^*$.

Since $V$ is affine, we can conclude that $Z^1 := Y \setminus V$ is of pure dimension one. For $i, j \in I$ let $Z^1_{ij} := Y_{ij} \setminus V$ and denote by $\mathfrak{D}^1_{ij}$ the set of irreducible components of $Z^1_{ij}$.

**6.2 Step.** *Choose a finite set $M^2 \subset \mathbb{C}[V]$ of homogeneous functions such that for all $i, j \in I$ and every $D \in \mathfrak{D}^1_{ij}$ there are $f_1, \ldots, f_r \in M^2 \cap \mathbb{C}[Y_{ij}]$ that generate the ideal*

$$I(D, i, j) := \{g \in \mathbb{C}[Y_{ij}];\ g|_D = 0\}.$$ □

Let $Z^2$ denote the union of all closures, taken with respect to $Y$, of the sets $N(V; f)$, where $f \in M^2$. For $i, j \in I$, let $Z^2_{ij} := Y_{ij} \cap Z^2$. Let $\mathfrak{D}^2_{ij}$ denote the set of irreducible components of $Z^2_{ij}$.

**6.3 Step.** *Choose a finite set $M^3 \subset \mathbb{C}[V]$ of homogeneous functions such that*

i) *for all $i, j \in I$ and every $D \in \mathfrak{D}^2_{ij}$ there are $f_1, \ldots, f_r \in M^3 \cap \mathbb{C}[Y_{ij}]$ such that $N(Y_{ij}; f_1, \ldots, f_r) = D$,*

ii) *for any two different $f, f' \in M^3$ the common components of $(f)$ and $(f')$ lie in $Z^1 \cup Z^2$.*

Note that the above step is the equivariant version of [10], Step 7. The following step corresponds to [10], Step 8:

**6.4 Step.** *Choose a finite set $M^4 \subset \mathbb{C}[V]$ of homogeneous functions such that*

i) *given a homogeneous $g \in \mathbb{C}(Y)^*$ and $p \in Y_{ij} \cap \mathrm{Def}(g)$ with $(g)_p \subset Z^1_{ij} \cup Z^2_{ij}$, there are functions $f_1, \ldots, f_r \in M^4 \cap \mathbb{C}[Y_{ij}]$ with $\alpha_0(g)_p = \alpha_1(f_1)_p + \ldots + \alpha_r(f_r)_p$, where $\alpha_i \in \mathbb{Z}_{>0}$,*

ii) *for any two different $f, f' \in M^3 \cup M^4$ the common components of $(f)$ and $(f')$ lie in $Z^1 \cup Z^2$.*

Now let $M^k$ be as in the above steps let $M$ denote the subgroup of $\mathbb{C}(Y)$ generated by $M^1 \cup \ldots \cup M^4$. For $i, j \in I$ let $\mathfrak{D}^3_{ij}$ denote the set of divisors obtained from $(f)$, $f \in M^3 \cup M^4$ by dropping the components belonging to $Z^1_{ij}$ or $Z^2_{ij}$. Set

$$\mathfrak{D}_{ij} := \mathfrak{D}^1_{ij} \cup \mathfrak{D}^2_{ij} \cup \mathfrak{D}^3_{ij}.$$

**6.5 Theorem.** *On each set $Y_{ij}$, the group $M \subset \mathbb{C}(Y_{ij})$ and the set $\mathfrak{D}_{ij}$ of divisors satisfy to SF0 – SF4.*

The proof of this Theorem is given in Section 7. For later use, we note some properties of the family of equivariant e.s.f. constructed above. For $i, j \in I$ and $p \in Y_{ij}$ let

$$U^{ij}_p := \{y \in Y_{ij}; \ \forall f \in S_p \ (f)_y \geq 0\}.$$

**6.6 Proposition.** *Let $M \subset \mathbb{C}(Y)$ and $\mathfrak{D}_{ij}$ be as above. Then we have for any two $i, j \in I$:*

i) *The set $Y_{ii} \setminus Y_{ij}$ is a union of irreducible divisors of $\mathfrak{D}_{ii}$.*

ii) *For every $p \in Y_{ij}$ the set $U^{ij}_p$ equals $U^{ii}_p$.*

iii) *For any $p \in Y_{ii}$ and $q \in Y_{jj}$ one obtains*

$$U^{ii}_p \cap U^{jj}_q = \bigcup_{y \in U^{ii}_p \cap U^{jj}_q} U^{ii}_y = \bigcup_{y \in U^{ii}_p \cap U^{jj}_q} U^{ij}_y = \bigcup_{y \in U^{ii}_p \cap U^{jj}_q} U^{jj}_y.$$

**Proof.** In order to obtain i), note that each $Y_{ij}$ is affine and hence $Y_{ii} \setminus Y_{jj} = Y_{ii} \setminus Y_{ij}$ is of pure codimension on in $Y_{ii}$. Moreover, $Y_{ii} \setminus Y_{jj}$ is a closed subset of $Y_{ii} \setminus V = Z^1_{ii}$. Consequently $Y_{ii} \setminus Y_{jj}$ is a union of irreducible components of $Z^1_{ii}$ which gives us the assertion.

We check ii). Clearly we have $U^{ij}_p \subset U^{ii}_p$. In order to show $U^{ii}_p \subset U^{ij}_p$, let $y \in Y_{ii} \setminus U^{ij}_p$. If $y \in Y_{ij}$, then, by definition of $U^{ij}_p$, there is a function $f \in S_p$ that is not regular at $y$. That means $y \notin U^{ii}_p$.

If $y \notin Y_{ij}$, then, by Property i) we have $y \in |D|$ with some $D \in \mathfrak{D}_{ii}$ that does not meet $U^{ij}_p$. Now, SF2 provides a function $f \in M$ with $(f)_p \leq 0$ that vanishes along $D$. Then $f^{-1}$ lies in $S_p$ and is not regular at $y$. Thus $y \notin U^{ii}_p$.

For the proof of assertion iii), set $U := U^{ii}_p \cap U^{jj}_q$. Using Lemma 5.1 and Property ii) we obtain

$$U = (U^{ii}_p \cap Y_{ij}) \cap (U^{jj}_q \cap Y_{ij})$$

$$= \bigcup_{y \in U^{ii}_p \cap Y_{ij}} U^{ii}_y \cap \bigcup_{y \in U^{jj}_q \cap Y_{ij}} U^{jj}_y$$

$$= \bigcup_{y \in U^{ii}_p \cap Y_{ij}} U^{ij}_y \cap \bigcup_{y \in U^{jj}_q \cap Y_{ij}} U^{ij}_y$$

$$= \bigcup_{y, y' \in U} (U^{ij}_y \cap U^{ij}_{y'})$$

$$= \bigcup_{y \in U} U^{ij}_y = \bigcup_{y \in U} U^{ii}_y = \bigcup_{y \in U} U^{jj}_y. \qquad \square$$

## 7 Proof of Steps 6.3 and 6.4

**Proof of Step 6.3.** Choose a total ordering of the set $S := \{(i,j,D);\ i,j \in I,\ D \in \mathfrak{D}_{ij}^2\}$. Suppose we have found an at most finite set $M' \subset \mathbb{C}[V]$ that satisfies 6.3 i) and ii) for all elements of $S$ that are strictly smaller than a given $(i,j,D)$. Let

$$V' := V \setminus Z_{ij}^2, \qquad A_1 := \bigcup_{f \in M'} N(V'; f), \qquad y_1 \in Y_{ij} \setminus D.$$

Applying Lemma 3.2 to $D$, $V'$, $A_1$ and $y_1$, we find an homogeneous $f_1 \in \mathbb{C}[Y_{ij}]$ with $f_1|_D = 0$ and $f_1(y_1) \neq 0$ such that $N(V'; f_1)$ contains no irreducible component of $A_1$. Set

$$W_1 := Y_{ij} \setminus (D \cup N(Y_{ij}; f_1)).$$

If $W_1$ equals $Y \setminus D$, then $M' \cup \{f\}$ satisfies 6.3 i) and ii) for all elements of $S$ less or equal to $(i,j,D)$. If there is a point $y_2 \in (Y_{ij} \setminus D) \setminus W_1 = N(Y_{ij}; f_1) \setminus D$, then set

$$A_2 := A_1 \cup N(V'; f_1)$$

and choose $f_2 \in \mathbb{C}[Y_{ij}]$ homogeneous such that $f_2|_D = 0$, $f_2(y_2) \neq 0$ and $N(V; f_2)$ contains no irreducible component of $A_2$. Let

$$W_2 := W_1 \cup (Y_{ij} \setminus (D \cup N(Y_{ij}; f_2))).$$

Note that $W_1$ is a proper subset of $W_2$. Defining $y_3$, $A_3$, $f_3$, $W_3$ etc., we obtain a strictly increasing sequence of open subsets $W_i$. Since $Y_{ij} \setminus D$ is a noetherian topological space, we have $W_l = Y_{ij} \setminus D$ for some $l$.

By construction, the set $M' \cup \{f_1, \ldots, f_l\}$ satisfies to Properties 6.3 i) and ii) for all elements of $S$ that are less or equal to $(i,j,D)$. So we can construct $M^3$ inductively. □

**Proof of Step 6.4.** At the beginning, we follow closely [10], p. 718. A family $(V_{ij})_{i,j \in I}$ of $\mathbb{C}^*$-stable open subsets $V_{ij} \subset Y_{ij}$ is called admissible, if $(V \setminus Z^2) \subset V_{ij}$ holds for any two $i,j \in I$ and there is a finite set $M' \subset \mathbb{C}[V]$ satisfying Property 6.4 ii) and as well Property 6.4 i) after replacing $Y_{ij}$ with $V_{ij}$.

Note that $V_{ij} := V \setminus Z^2$ and $M' := \{1\}$ defines an admissible family. Since $Y$ is a noetherian topological space, we can choose an admissible family $(V_{ij})_{i,j \in I}$ such that each $V_{ij}$ is maximal. We show that in this case $V_{ij} = Y_{ij}$ holds for any two $i,j \in I$.

Let $M' \subset \mathbb{C}[V]$ be a set of functions making the above family $(V_{ij})_{i,j \in I}$ admissible. Assume that we have $V_{ij} \neq Y_{ij}$ for some fixed $i,j \in I$. Set $\mathfrak{C} := \mathfrak{D}_{ij}^1 \cup \mathfrak{D}_{ij}^2$ and for $y \in Y_{ij}$ let $\mathfrak{C}^y := \{D \in \mathfrak{C};\ y \in D\}$. Moreover, let

$$p \in A := Y_{ij} \setminus V_{ij} = (Z_{ij}^1 \cup Z_{ij}^2) \setminus V_{ij}$$

be a point for that $\mathfrak{C}^p$ has a minimal number of elements. Then, by removing the elements of $\mathfrak{C} \setminus \mathfrak{C}^p$ from $Y_{ij}$, we obtain a $\mathbb{C}^*$-invariant open neighbourhood $U \subset Y_{ij}$ of the point $p$ such that $\mathfrak{C}^y = \mathfrak{C}^p$ holds for every $y \in A \cap U$.

Let W (resp. $W^+$) denote the free abelian group (resp. semigroup) of Weil divisors (resp. effective Weil divisors) generated by the elements of $\mathfrak{C}$. For $y \in A \cap U$ let $\mathrm{P}(y)$ (resp. $\mathrm{P}^+(y)$) denote the set of all $D \in W$ (resp. $D \in W^+$) such that the components of $D$ belong to $\mathfrak{C}^y$ and there is a homogeneous $g \in \mathbb{C}^*(Y)$ with $D = (g)_y$.

Since $P(y)$ is a subgroup of $W$ and $P^+(y) = W^+ \cap P(y)$, we see by Gordon's Lemma that each semigroup $P^+(y)$ is finitely generated. In particular, for each $y \in U \cap A$ there is a $\mathbb{C}^*$-stable neighbourhood $U^y \subset U$ of $y$ such that $P^+(y) \subset P^+(u)$ for all $u \in U^y \cap A$.

By the elementary Lemma 7.1 stated below, we can assume that the semigroup $P^+(p)$ is maximal among the $P^+(y)$, $y \in U \cap A$. Consequently, by shrinking $U$ appropriately, we achieve that $P^+(y) = P^+(p)$ holds for all $y \in U \cap A$.

After shrinking $U$ again, we find homogeneous functions $h_1, \ldots, h_r \in \mathcal{O}(U)$ such that the divisors $(h_k)_p$ generate $P^+(p)$. In order to ensure Property 6.4 ii) for these functions, we are going to modify them slightly.

Suppose we found homogeneous functions $\widetilde{h}_1, \ldots, \widetilde{h}_{s-1} \in \mathbb{C}[Y_{ij}]$ such that $(\widetilde{h}_k)_p = \alpha_k (h_k)_p$ holds with some $\alpha_k \in \mathbb{Z}_{>0}$ and 6.4 ii) is satisfied for $M^3 \cup M' \cup \{\widetilde{h}_1, \ldots, \widetilde{h}_{s-1}\}$, i.e. for any two different functions $f$ and $f'$ of this set, the common components of the divisors $(f)$ and $(f')$ are contained in $Z^1 \cup Z^2$. Set $V' := V \setminus Z^2$ and

$$B := \bigcup_{f \in M^3} \overline{N(V'; f)} \cup \bigcup_{f \in M'} \overline{N(V'; f)} \cup \bigcup_{k=1}^{s-1} \overline{N(V'; \widetilde{h}_k)}.$$

Here the closures are taken with respect to $Y_{ij}$. Since $(h_s)_p$ has only components in $\mathfrak{C}^p$, no irreducible component $B_0$ of $B$ with $p \in B_0$ occurs in $(h_s)$. Let $B_n$, $n = 1, \ldots, q$, denote the irreducible components of $B$ with $p \notin B_n$. According to Lemma 3.2, choose homogeneous functions $b_n \in \mathbb{C}[Y_{ij}]$ such that

$$b_n(p) \neq 0, \qquad b_{n_1}|_{B_{n_2}} = 0 \Leftrightarrow n_1 = n_2.$$

Modify $h_s$ step by step as follows

$$_1h_s := \frac{h_s^{\mathrm{ord}_{B_1}(b_1)}}{b_1^{\mathrm{ord}_{B_1}(h_s)}}, \quad \ldots, \quad _qh_s := \frac{_{q-1}h_s^{\mathrm{ord}_{B_q}(b_q)}}{b_q^{\mathrm{ord}_{B_q}(_{q-1}h_s)}} =: h'_s.$$

note that $(h'_s)_p = \alpha_s (h_s)_p$ holds for some $\alpha_s \in \mathbb{Z}_{>0}$. By construction, the divisor $(h'_s)$ has no components that are contained in $B$. We make $h'_s$ to a regular function on $Y_{ij}$ by means of the following modification:

Let $A_l \subset Y_{ij}$ denote the components of the pole divisor of $h'_s$, viewed as a rational function on $Y_{ij}$. Lemma 3.2 provides homogeneous functions $a_l \in \mathbb{C}[Y_{ij}]$ such that

$$a_l(p) \neq 0, \qquad a_l|_{B_n} \neq 0, n = 1, \ldots, q, \qquad \mathrm{ord}_{A_l}(a_l) > -\mathrm{ord}_{A_l}(h'_s).$$

Set $\widetilde{h}_s := h'_s \prod_l a_l$. Then $\widetilde{h}_s$ is regular on $Y_{ij}$, does not vanish along any irreducible component of $B$ and satisfies $(\widetilde{h}_s)_p = \alpha_s (h_s)_p$.

Proceeding this way we arrive finally at $s = r$. Since $(\widetilde{h}_k)_p = \alpha_k (h_k)_p$ holds with certain $\alpha_k \in \mathbb{Z}_{>0}$, for every $D \in P^+(y)$, $y \in U \cap A$, some positive multiple of $D$ is a nonnegative integral linear combination of the $(h_k)_p$.

Thus, replacing $V_{ij}$ with $V_{ij} \cup U$ and $M'$ with $M' \cup \{\widetilde{h}_1, \ldots, \widetilde{h}_r\}$ yields an admissible family $(\widetilde{V}_{kl})_{k,l \in I}$ such that $\widetilde{V}_{ij}$ properly contains $V_{ij}$. This contradicts the maximality property of the admissible family we started with. □

**7.1 Lemma.** *Let $N := \mathbb{Z}^n$ and let $N^+$ denote the positive quadrant in $N$. If $N_i$, $i \in I$, is a family of sublattices of $N$, then the family $(N_i \cap N^+)$, $i \in I$, contains maximal elements with respect to inclusion.* □

# 8  Proof of Theorem 6.5

**8.1 Lemma.** *For any two $i, j \in I$, there is a function $h \in \mathbb{C}[Y_{ij}] \cap M$ with $N(Y_{ij}; h) = Z_{ij}^1 \cup Z_{ij}^2$.*

**Proof.** For each $D \in \mathfrak{D}_{ij}^2$ choose an $h_D \in M^2$ that vanishes along $D$. Note that $h_D$ is regular on $Y_{ij} \setminus Z_{ij}^1 = Y_{ij} \cap V$. Set

$$h' := \prod_{D \in \mathfrak{D}_{ij}^2} h_D.$$

Then $h'$ is regular on $Y_{ij} \setminus Z_{ij}^1$ and the components of $(h')|_{Y_{ij}}$ lie in $\mathfrak{D}_{ij}^1 \cup \mathfrak{D}_{ij}^2$. Thus Step 6.2, provides a function $h'' \in M^2 \cap \mathbb{C}[Y_{ij}]$ such that $h := h'h''$ is as wanted. □

**Proof of Theorem 6.5.** Let $i, j \in I$ be fixed. Since the functions $f \in M^3 \cup M^4$ are regular on $V$, every divisor of $\mathfrak{D}_{ij}$ is effective. Moreover, SF0 and SF1 are valid by definition.

Let us check SF2. So, let $D \in \mathfrak{D}_{ij}$ and let $y \in Y_{ij} \setminus |D|$. If $D \in \mathfrak{D}_{ij}^1$, then Step 6.2 yields an $f \in M^2 \cap \mathbb{C}[Y_{ij}]$ with $f|_D = 0$ and $f(y) \neq 0$. For $D \in \mathfrak{D}_{ij}^2$, the desired function $f$ is provided by Step 6.3.

Now, let $D \in \mathfrak{D}_{ij}^3$. Then there is an $f_1 \in M^3 \cup M^4$ with $(f_1) = D + \widetilde{D}$ with a divisor $\widetilde{D}$ supported in $Z_{ij}^1 \cup Z_{ij}^2$. Let $h \in \mathbb{C}[Y_{ij}]$ be as in Lemma 8.1. Then, for a suitable power $s > 0$, we have $h_q^s \geq (f_1)_q$. Thus $h^s f_1^{-1}$ is regular at $q$.

By Step 6.4, we find an $f_2 \in M$ with $(f_2)_q = n(h^s f_1^{-1})_q$ for some positive integer $n$. Consider $f_3 := f_1^n f_2 h^{-ns} \in M$. Then $(f_3)_q = 0$. Moreover,

$$(f_3)_p = nD + \sum_{D' \in \mathfrak{D}_{ij}^1 \cup \mathfrak{D}_{ij}^2} n_{D'} D' + \sum_{D'' \in \mathfrak{D}_{ij}^3} n_{D''} D''$$

with all $n_{D''} \geq 0$ and $n_{D'} = 0$ if $q \in D'$. Now, let $\mathfrak{D}_0$ be the set of divisors $D'$ in $\mathfrak{D}_{ij}^1 \cup \mathfrak{D}_{ij}^2$ with $n_{D'} = 0$ in the above representation.

Since $D' \in \mathfrak{D}_0$ implies $q \notin D'$, the first two cases yield for every $D' \in \mathfrak{D}_0$ a function $f_{D'} \in M$ with $(f)_p \geq D'$ and $(f)_q = 0$. Then the desired function is

$$f := f_3 \prod_{D' \in \mathfrak{D}_0} f_{D'}.$$

Next we verify SF3: Let $p \in Y_{ij}$ and let $g \in \mathbb{C}(Y_{ij})^*$ be homogeneous such that $(g)_p = \sum_{D \in \mathfrak{D}_{ij}} n_D D_p$. We can write this divisor as

$$(g)_p = \left( \sum_{D' \in \mathfrak{D}_{ij}^1 \cup \mathfrak{D}_{ij}^2} n_{D'} D' \right) + \left( \sum_{D'' \in \mathfrak{D}_{ij}^3} n_{D''} (D'')_p \right).$$

For every $D'' \in \mathfrak{D}_{ij}^3$ there is by definition a function $h_{D''} \in M^3 \cup M^4$ such that the components of the divisor $(h_{D''}) - D''$ belong to $\mathfrak{D}_{ij}^1 \cup \mathfrak{D}_{ij}^2$. Set

$$h_1 := \prod_{D'' \in \mathfrak{D}_{ij}^3} h_{D''}^{-n_{D''}}.$$

By construction, all components of the divisor $(gh_1)_p$ lie in $\mathfrak{D}_{ij}^1 \cup \mathfrak{D}_{ij}^2$. Moreover, choose a function $h_2 \in \mathbb{C}[Y_{ij}] \cap M$ as in Lemma 8.1.

Then, for a suitable power $s > 0$, the function $gh_1h_2^s$ is regular at $p$ and the components of $(gh_1h_2^s)_p$ lie in $\mathfrak{D}_{ij}^1 \cup \mathfrak{D}_{ij}^2$. Thus Step 6.4 yields an $f' \in M$ such that $\alpha(gh_1h_2^s)_p = (f')_p$ for some $\alpha \in \mathbb{Z}_{>0}$. Then $f := f'h_1^{-1}h_2^{-s}$ is as wanted.

Finally, we come to the verification of Property SF4. So, let $p \in Y_{ij}$. In a first step we show that

$$(*) \qquad W_p := Y_{ij} \setminus \left( \bigcup_{D \in \mathfrak{D}_{ij} \setminus \mathfrak{D}_{ij}^p} |D| \right) = U_p.$$

In order to verify the inclusion "$\subset$", let $y \in W_p$. By definition of $W_p$, we have $\mathfrak{D}_{ij}^y \subset \mathfrak{D}_{ij}^p$ and hence $\sigma_y \subset \sigma_p$. This implies $S_p \subset S_y$, which in turn means $U_y \subset U_p$. In particular it follows $y \in U_p$.

The inclusion "$\supset$" is shown as follows. Consider a point $y \in Y_{ij} \setminus W_p$. Then $y \in |D|$ for some $D \in \mathfrak{D}_{ij} \setminus \mathfrak{D}_{ij}^p$. We have to treat the following cases:

*Case 1:* $D \in \mathfrak{D}_{ij}^1$. Then Step 6.2 provides a function $f \in \mathbb{C}[Y_{ij}]$ with $f|_D = 0$ and $f(p) \neq 0$. Then $f^{-1} \in S_p$ and $f^{-1} \notin S_y$, which yields $y \notin U_p$.

*Case 2:* $D \in \mathfrak{D}_{ij}^2$. Here we find by Step 6.3 a function $f \in \mathbb{C}[Y_{ij}]$ with $f|_D = 0$ and $f(p) \neq 0$. As in Case 1, this leads to $y \notin U_p$.

*Case 3:* $D \in \mathfrak{D}_{ij}^3$. By definition, there is a function $f \in M^3 \cup M^4$ such that

$$(f)|_{Y_{ij}} = D + \sum_{D' \in \mathfrak{D}_{ij}^1 \cup \mathfrak{D}_{ij}^2} n_{D'} D'.$$

Let $h \in M \cap \mathbb{C}[Y_{ij}]$ be as in Lemma 8.1. Thus, for an appropriate power $s > 0$, the function $s^n f^{-1} \in M$ is regular at $p$ but not at $y$. This means $y \notin U_p$.

Now, having verified the representation $(*)$, we can show that $U_p$ is affine. By the above Cases 1 and 2, we find for each $D \in \mathfrak{D}_{ij}^1 \cup \mathfrak{D}_{ij}^2 \setminus \mathfrak{D}_{ij}^p$ a function $f_D \in \mathbb{C}[Y_{ij}]$ that vanishes along $D$ and does not vanish at $p$. Hence we obtain

$$U_p \subset \widetilde{U}_p := \bigcap_{D \in (\mathfrak{D}_{ij}^1 \cap \mathfrak{D}_{ij}^2) \setminus \mathfrak{D}_{ij}^p} (Y_{ij})_{f_D}$$

For each $D \in \mathfrak{D}_{ij}^3 \setminus \mathfrak{D}_{ij}^p$ choose an $f'_D \in M^3 \cup M^4$ with $(f'_D) = D + \widetilde{D}$ with some divisor $\widetilde{D}$ supported in $Z_{ij}^1 \cup Z_{ij}^2$. Moreover, let $h \in \mathbb{C}[Y_{ij}]$ be as in in Lemma 8.1.

For a suitable power $s > 0$, the divisor $(f'_D h^{-s})_p$ has only negative components. So, Step 6.4 gives us a function $g \in M \cap \mathbb{C}[Y_{ij}]$ with $(g)_p = -(f'_D h^{-s})_p$. Then $f_D := f'_D g$ is regular on $\widetilde{U}_p$, vanishes along $D$ and does not vanish at $p$. Now, $(*)$ implies

$$U_p \subset \bigcap_{D \in \mathfrak{D}_{ij}^3 \setminus \mathfrak{D}_{ij}^p} (\widetilde{U}_p)_{f_D} \subset W_p = U_p.$$

In particular, $U_p$ is affine. Moreover, since $M^1$ contains generators of $\mathbb{C}[Y_{ij}]$, we see that $\mathbb{C}[U_p]$ is generated by

$$(M^1 \cap \mathbb{C}[Y_{ij}]) \cup \{f_D^{-1}; \ D \in \mathfrak{D}_{ij} \setminus \mathfrak{D}_{ij}^p\} \subset \mathbb{C}\{S_p\}. \qquad \square$$

# 9    The Main Result

**9.1 Theorem.** *Every normal complex algebraic $\mathbb{C}^*$-variety $Y$ admits an equivariant closed embedding into a toric prevariety $X$ on which $\mathbb{C}^*$ acts as a one-parameter-subgroup of the big torus $T \subset X$. If $Y$ is $\mathbb{Q}$-factorial, then $X$ may be chosen to be simplicial and of affine intersection.*

Our proof of the second part of Theorem 9.1 is based on the fact that in the construction there is a certain freedom for glueing the embeddings obtained by a family of equivariant e.s.f. in the $\mathbb{Q}$-factorial case. We demonstrate this phenomenon in the following

**9.2 Example.** We construct a certain e.s.f. on $Y := \mathbb{C}^2$. Consider the multiplicative subgroup $M$ of $\mathbb{C}(Y)$ generated by the functions

$$f_1(z,w) := z-1, \quad f_2(z,w) := z+1, \quad f_3(z,w) := w+z^2-1, \quad f_4(z,w) := w.$$

Note that $M \cong \mathbb{Z}^4$. Set $D_i := (f_i) = N(Y; f_i)$ and let $\mathfrak{D} := \{D_1, \ldots, D_4\}$. Then we have the following real picture for the divisors $D_i$.

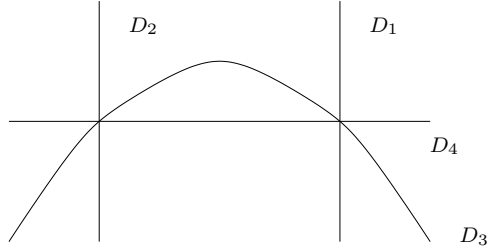

By definition, the divisors $D_i$ are irreducible and hence two distinct $D_i$ and $D_j$ don't have common components. Moreover, SF0, SF1 and SF2 are obviously satisfied. SF3 follows from the fact $I(Y; D_i) = (f_i)$. In order to check SF4, we have to deal with the following four cases:

*Case 1:* $p \in Y \setminus \bigcup_i D_i$. Then $U_p = Y \setminus \bigcup_i D_i$ and for $f := \prod_i f_i$ one has $\mathcal{O}(U_p) = \mathbb{C}[z,w]_f$.

*Case 2:* $p \in D_i \setminus \bigcup_{j \ne i} D_j$. Then $U_p = Y \setminus \bigcup_{j \ne i} D_j$ holds and setting $f := \prod_{j \ne i} f_i$ one has $\mathcal{O}(U_p) = \mathbb{C}[z,w]_f$.

*Case 3:* $p = (-1, 0)$. Then we have $U_p = Y \setminus D_1$ and $\mathcal{O}(U_p) = \mathbb{C}[z,w]_{f_1}$.

*Case 4:* $p = (1, 0)$. Then we have $U_p = Y \setminus D_2$ and $\mathcal{O}(U_p) = \mathbb{C}[z,w]_{f_2}$.

So we see that $M$ is an e.s.f. on $Y$. The fan $\Delta(M, Y)$ has precisely two maximal cones, namely $\sigma_{(1,0)}$ and $\sigma_{(-1,0)}$. On the other hand, $M$ defines also e.s.f.'s on the sets $Y_{11} := Y \setminus D_2$, $Y_{22} := Y \setminus D_1$ and $Y_{12} := Y_1 \cap Y_2$. We have

$$\Delta(M, Y_{11}) = \mathfrak{F}(\sigma_{(1,0)}), \quad \Delta(M, Y_{22}) = \mathfrak{F}(\sigma_{(1,0)}), \quad \Delta(M, Y_{12}) = \{\{0\}, \sigma_{(0,0)}, \sigma_{(0,1)}\}.$$

Here $\mathfrak{F}(\sigma)$ denotes the fan of faces of a cone $\sigma$. Since the fan $\Delta(M, Y_{12})$ has two maximal cones, the glueing of the maps $\varphi_{M, Y_{ii}}$ along $Y_{12}$ embeds $Y$ into a toric prevariety that is not of affine intersection. ◇

For the proof of Theorem 9.1, we recall the following notions from the theory of toric varieties. Let $\Delta$ be a fan in some lattice and let $\tau$ be any cone of $\Delta$. We denote by $x_\tau$ the distinguished point of $X_\Delta$ associated to $\tau$ (see [6], p. 28) and by $V_\tau$ the closure of the orbit $T \cdot x_\tau \subset X_\Delta$ (see [6], Section 3.1).

If $\Delta = \Delta(M, Y)$ with some equivariant s.f. on an $H$-variety $Y$ and $\mathfrak{D}$ is a set of divisors on $Y$ satisfying SF0) – SF3), then we define for each $D \in \mathfrak{D}$ a one-dimensional cone $\varrho_D := \mathbb{R}_{\geq 0}\mathrm{val}_D$. Recall that by Lemma 4.4 the $\varrho_D$ are precisely the extremal rays of the fan $\Delta(M, Y)$.

**9.3 Lemma.** *Let $\mathfrak{D}$ be a set of divisors that satisfies SF0) – SF3) for some equivariant s.f. $M \subset \mathbb{C}(Y)$. Then $D \in \mathfrak{D}$ implies $\varphi_{M,Y}^{-1}(V_{\varrho_D}) = |D|$.*

**Proof.** Let $y \in |D|$ such that $\mathfrak{D}^y = \{D\}$ and choose $y_0 \in U_y$ with $\mathfrak{D}^{y_0} = \emptyset$. Let $T$ denote the big torus of $X_{\Delta(M,Y)}$. Since one has $X_{\varrho_D} = T \cdot x_0 \cup T \cdot x_{\varrho_D}$ (see [6], p. 54), Proposition 5.5 i) yields

$$\varphi_{M,Y}^{-1}(T \cdot x_0) = U_{y_0}, \qquad \varphi_{M,Y}^{-1}(T \cdot x_{\varrho_D}) = U_y \setminus U_{y_0} = U_y \cap |D|.$$

This implies $\varphi_{M,Y}(U_y \cap |D|) \subset T \cdot x_{\varrho_D}$. Hence we have $\varphi_{M,Y}(|D|) \subset V_{\varrho_D}$ which means $|D| \subset \varphi_{M,Y}^{-1}(V_{\varrho_D})$. To obtain $\varphi_{M,Y}^{-1}(V_{\varrho_D}) \subset |D|$, consider a point $y \in Y \setminus |D|$. Then $D \notin \mathfrak{D}^y$ and hence $\varrho_D \not\prec \sigma_y$, which implies $X_{\sigma_y} \cap V_{\varrho_D} = \emptyset$. So, we obtain $f(M,Y)(y) \notin V_{\varrho_D}$. □

**9.4 Lemma.** *Suppose that $M \subset \mathbb{C}(Y)$ defines an equivariant weak s.f. on the $H$-variety $Y$ and also on an $H$-stable open subset $V \subset Y$. Then*

i) $\Delta(M, V)$ *is a subfan of* $\Delta(M, Y)$.

ii) $\varphi_{M,Y}|_V = \varphi_{M,V}$.

**Proof.** Statement i) follows from recalling the definitions: Let $\sigma_p \in \Delta(M, V)$. Then $\sigma_p = \mathrm{cone}(S_p)^\vee$, where $S_p = \{f \in M;\ (f)_p \geq 0\}$, i.e., $\sigma_p \in \Delta(M, Y)$.

To show ii), let $T$ denote the acting torus of the toric variety $X_{\Delta(M,Y)}$ and let $p \in V$ be a point such that each $f \in M$ is regular in $p$. Then Proposition 5.5 i) yields

$$\begin{aligned} \varphi_{M,Y}^{-1}(T) &= U_p = \{y \in Y;\ \forall f \in M\ (f)_y = 0\}, \\ \varphi_{M,V}^{-1}(T) &= \{y \in V;\ \forall f \in M\ (f)_y = 0\} = U_p \cap V. \end{aligned}$$

It suffices to show that $\varphi_{M,Y}$ and $\varphi_{M,V}$ coincide on the (nonempty) open set $U_p \cap V$. But this is a consequence of the fact that the inclusions $M \subset \mathcal{O}(U_p) \subset \mathcal{O}(U_p \cap V)$ define a commutative diagram

$$\begin{array}{ccc} & \mathcal{O}(T \cdot x_0) & \\ \swarrow & & \searrow \\ \mathcal{O}(U_p) & \xrightarrow[f \mapsto f|_{U_p \cap V}]{} & \mathcal{O}(U_p \cap V) \end{array}. \quad \square$$

**Proof of Theorem 9.1.** According to Sumihiro's Theorem, choose a cover of $Y$ by finitely many open affine $\mathbb{C}^*$-stable subsets $Y_i$, $i \in I$. Let $M \subset \mathbb{C}(Y)$ be a family that defines an equivariant e.s.f. on each $Y_{ij} := Y_i \cap Y_j$ constructed by the procedure of Section 6.

We define an affine system of fans in $N := \mathrm{Hom}(M, \mathbb{Z})$ as follows: Let $K \subset Y \times I$ be a finite set such that for each fixed $i \in I$ the cones $\sigma_p$, $(p, i) \in K$, are just the maximal cones of $\Delta(M, Y_{ii})$. For two elements $k = (p, i)$ and $k' = (q, j)$ of $K$ set

$$\Delta_{kk'} := \Delta_{pq}^{ij} := \bigcup_{y \in U_p^{ii} \cap U_q^{jj}} \mathfrak{F}(\sigma_y).$$

Then $\mathcal{S} := (\Delta_{kk'})_{k,k' \in K}$ is an affine system of fans in the lattice $N$. Note that $\Delta_{pq}^{ii}$ is the fan of faces of $\sigma_p$. Now, for every index $k = (p,i) \in K$, set

$$U_k := U_p^{ii}, \qquad \varphi_k := \varphi_{M,p}^{ii} \colon U_k \to X_{\sigma_p}.$$

Note that $Y$ is covered by the $U_k$. Let $T$ denote the acting torus of the $X_{\sigma_p}$ and recall that there is a one-parameter-subgroup $\lambda \colon \mathbb{C}^* \to T$ making all $\varphi_k$ equivariant. Proposition 6.6 and Lemma 5.4 i) imply that for any two $k, k'$ in $K$ there is a commutative diagram

$$\begin{array}{ccccc}
U_k & \supset & U_k \cap U_{k'} & \subset & U_{k'} \\
\varphi_k \downarrow & & \downarrow & & \downarrow \varphi_{k'} \\
X_{\Delta_{kk}} & \supset & X_{\Delta_{kk'}} & \subset & X_{\Delta_{k'k'}}
\end{array}.$$

Consequently, the maps $\varphi_k$ glue together to a regular map $\varphi \colon Y \to X_{\mathcal{S}}$ that is equivariant with respect to the $\mathbb{C}^*$-action on $X_{\mathcal{S}}$ defined by $\lambda$. Moreover, by Proposition 5.5 i) and Proposition 6.6 iii), we have for any two $k, k'$ in $K$ the identity

$$(\varphi_k)^{-1}(X_{\Delta_{kk'}}) = U_k \cap U_{k'} = (\varphi_{k'})^{-1}(X_{\Delta_{kk'}})$$

In particular, it follows $\varphi^{-1}(X_{[\sigma_p,k]}) = U_k$ for any $k = (p,i) \in K$. Since, by Lemma 5.4 ii), each $\varphi_k$ is a closed embedding, also $\varphi$ is a closed embedding. Hence the first part of the Main Theorem is proved.

Now assume that $Y$ is $\mathbb{Q}$-factorial. Then Lemma 4.6 yields that $X_{\mathcal{S}}$ is simplicial. However, $X_{\mathcal{S}}$ need in general not be of affine intersection. We consider the reduction $g \colon X_{\mathcal{S}} \to X_{\widetilde{\mathcal{S}}}$ of Proposition 1.1 and set $\eta := g \circ \varphi$.

Then $\eta$ also is equivariant with respect to the $\mathbb{C}^*$-action on $X_{\widetilde{\mathcal{S}}}$ defined by the one-parameter-subgroup $\lambda$. We claim that $\eta$ still is an embedding. In order to verify this, we only have to check that $\eta^{-1}(X_{[\widetilde{\sigma}_{k'k'},k']}) = U_{k'}$ holds for every $k' \in K$. By Proposition 1.1 ii) one has

$$\eta^{-1}(X_{[\widetilde{\sigma}_{k'k'},k']}) = \varphi^{-1}(g^{-1}(X_{[\widetilde{\sigma}_{k'k'},k']}))$$
$$= \varphi^{-1}\left(\bigcup_{k \in K} X_{[\widetilde{\sigma}_{kk'},k]}\right).$$

For fixed $k, k'$ we have

$$\varphi^{-1}(X_{[\widetilde{\sigma}_{kk'},k]}) = \bigcup_{k'' \in K} \varphi|_{U_{k''}}^{-1}(X_{[\widetilde{\sigma}_{kk'},k]})$$
$$= \bigcup_{k'' \in K} \varphi|_{U_{k''}}^{-1}(X_{[\widetilde{\sigma}_{kk'},k]} \cap X_{[\sigma_{k''k''},k'']})$$
$$= \bigcup_{k'' \in K} \varphi|_{U_{k''}}^{-1}\left(\bigcup_{\tau \in \Delta_{kk''};\, \tau \prec \widetilde{\sigma}_{kk'}} X_{[\tau,k'']}\right)$$
$$\subset \bigcup_{k'' \in K} \varphi|_{U_{k''}}^{-1}\left(\bigcup_{\tau \prec \widetilde{\sigma}_{kk''};\, \tau \prec \widetilde{\sigma}_{kk'}} X_{[\tau,k'']}\right)$$
$$\subset \bigcup_{k'' \in K} \varphi|_{U_{k''}}^{-1}\left(\bigcup_{\tau \prec \widetilde{\sigma}_{k'k''}} X_{[\tau,k'']}\right)$$

$$= \bigcup_{k'' \in K} \varphi|_{U_{k''}}^{-1}(X_{[\widetilde{\sigma}_{k'k''}, k'']})$$

$$= \bigcup_{k'' \in K} \varphi_{k''}^{-1}(X_{\widetilde{\sigma}_{k'k''}}).$$

Hence, to conclude the proof, it suffices to check for any two $k'' = (p, i)$ and $k' = (q, j)$ in $K$ the following condition

(∗) $$\varphi_{k''}^{-1}(X_{\widetilde{\sigma}_{k'k''}}) \subset U_{k'}$$

Recall from the definitions and Section 1 that $\widetilde{\sigma}_{k'k''} = \mathrm{conv}(\sigma_y; y \in U_{k'} \cap U_{k''})$. In order to verify (∗) consider a point $z \in U_{k''}$ with $\varphi_{k''}(z) \in X_{\widetilde{\sigma}_{k'k''}}$. We distinguish the following two cases:

*Case 1.* $z \in Y_{ij}$. By Lemma 9.4, we have $\varphi_{k''}(z) = \varphi_{M,Y_{ii}}(z) = \varphi_{M,Y_{ij}}(z)$. Applying Lemma 9.4 again, one obtains $\varphi_{M,Y_{jj}}(z) = \varphi_{M,Y_{ij}}(z) \in X_{\widetilde{\sigma}_{k'k''}} \subset X_{\sigma_q}$. Hence Proposition 5.5 i), applied to $\varphi_{M,Y_{jj}}$ yields $z \in U_{k'}$, which we wanted to obtain.

*Case 2.* $z \notin Y_{ij}$. By Proposition 6.6, we find an irreducible divisor $D \in \mathfrak{D}_{ii}$ such that $z \in |D| \subset Y_{ii} \setminus Y_{jj}$. Then $D \notin \mathfrak{D}_{ii}^y$ for every $y \in U_{k''} \cap U_{k'}$, i.e., $\varrho_D$ is not contained in any cone $\sigma_y$, with $y \in U_{k''} \cap U_{k'}$. Moreover, Lemma 9.3 yields $\varphi_{k''}(z) \in V_{\varrho_D} \subset X_{\Delta_{k''k''}}$, i.e.,

$$\varphi_{k''}(z) \in \bigcup_{\tau \in \Delta_{k''k''};\ \varrho_D \prec \tau} T \cdot x_\tau,$$

Since $\varrho_D$ is not a ray of $\widetilde{\sigma}_{k'k''}$, it follows that $V_{\varrho_D} \cap X_{\widetilde{\sigma}_{k'k''}} = \emptyset$. In particular, we have $\varphi_{k''}(z) \notin X_{\widetilde{\sigma}_{k'k''}}$, which is a contradiction. □

We conclude this section with two questions arising naturally from our Main Result: Is there an analogous equivariant embedding result for actions of tori of higher dimension or even for actions of reductive groups? Are there any nice necessary and sufficient conditions for existence of equivariant embeddings into toric varieties or into toric prevarieties of affine intersection?

## 10   Global Resolutions of Coherent Sheaves

Let $Y$ be a normal complex algebraic variety endowed with an effective regular action of a (possibly trivial) algebraic torus $H$. Moreover assume that $\mathcal{F}$ is a coherent $H$-sheaf on $Y$. For several applications it is useful to know, if $\mathcal{F}$ can be globally resolved by nice sheaves.

For trivial $H$, existence of locally free resolutions in the case of $Y$ being locally factorial is due to Kleiman (see e.g. [4]). For smooth $Y$, it is indicated in [8], Exercise III.6.9, how to obtain finite locally free resolutions. In [5], certain finite equivariant resolutions are obtained on toric varieties $Y$.

In this section we consider the following situation: Assume that there is an $H$-equivariant embedding $\imath: Y \to X$ into a toric prevariety $X$ of affine intersection, where $H$ acts on $X$ by means of a monomorphism $\varphi$ from $H$ to the acting torus $T$ of $X$. We prove

**10.1 Theorem.** *There exist $T$-stable Weil Divisors $D_1, \ldots, D_r$ on $X$ such that setting $\mathcal{E} := \bigoplus \imath^* \mathcal{O}_X(D_i)$ one obtains an exact sequence $\mathcal{E} \longrightarrow \mathcal{F} \longrightarrow 0$ of $H$-sheaves. If $X$ is smooth, then $\mathcal{E}$ is locally free.*

Here we mean by a Weil divisor on a prevariety a finite integral linear combination of closed irreducible one-codimensional subspaces. Since vanishing order is defined locally, we obtain as in the case of varieties a sheaf $\mathcal{O}(D)$ associated to each Weil divisor $D$ on a normal prevariety.

Note that the above result implies in particular the resolution result of Kleiman: If $Y$ is locally factorial, then the proof of [10], Lemma 3.11, and our construction (Section 9 and 1.1 iv)) yield that $Y$ can be embedded into a smooth toric prevariety of affine intersection. Our proof of Theorem 10.1 consists basically in adapting the arguments of [5], Section 1 to the situation of toric prevarieties of affine intersection.

First recall that by [2], Proposition 8.7, there exist an open toric subvariety $U$ of some $\mathbb{C}^n$ and a toric morphism $p\colon U \to X$ such the kernel $H_0$ of the homomorphism $\pi\colon (\mathbb{C}^*)^n \to T$ associated to $p$ is a torus and $p$ is a good prequotient (see [2], Section 6) for the action of $H_0$ on $U$.

We need some more detailed information about this map $p$. So, let $X$ arise from an affine system $\mathcal{S} = (\Delta_{ij})_{i,j \in I}$ of fans in a lattice $N$. Since $X$ is of affine intersection, each $\Delta_{ij}$ is the fan of faces of a single cone $\sigma_{ij}$. Denote by $R$ the set of equivalence classes $[\varrho, i] \in \Omega(\mathcal{S})$ with $\varrho$ one-dimensional. Let $N' := \mathbb{Z}^R$ and $\widetilde{N} := N \oplus N'$. Obtain lattice homomorphisms

$$P'\colon N' \to N, \quad P'(e_{[\varrho,i]}) := v_\varrho, \qquad P := \mathrm{id}_N + P'\colon \widetilde{N} \to N,$$

where $v_\varrho$ is the primitive lattice vector contained in $\varrho \in \bigcup \Delta_{ii}^{(1)}$ and $e_{[\varrho,i]}$ is the canonical base vector of $N'$ corresponding to $[\varrho, i] \in \Omega(\mathcal{S})$. For every $i \in I$ and $\sigma \in \Delta_{ii}$ define a strictly convex cone $\widetilde{\sigma}(i)$ in $\widetilde{N}$ by setting

$$\widetilde{\sigma}(i) := \mathrm{cone}(e_{[\varrho,i]};\ \varrho \in \sigma^{(1)})$$

Then the cones $\widetilde{\sigma_{ii}}(i)$, $i \in I$, are the maximal cones of a fan $\widetilde{\Delta}$ in $\widetilde{N}$ and $U$ is the toric variety defined by $\widetilde{\Delta}$. Moreover, the toric morphism $p$ arises from $P$ and we have

$$p^{-1}(X_{[\sigma,i]}) = X_{\widetilde{\sigma}(i)}$$

for each $\sigma \in \Delta_{ii}$, $i \in I$. Let $\widetilde{M} := \mathrm{Hom}(\widetilde{N}, \mathbb{Z})$. For $\widetilde{u} \in \widetilde{M}$ denote by $\chi^{\widetilde{u}}$ the corresponding character of $\widetilde{T} := (\mathbb{C}^*)^n$. Let $\mathcal{O}_U \otimes \mathbb{C}(\widetilde{u})$ denote the sheaf $\mathcal{O}_U \otimes \mathbb{C}$ endowed with the $\widetilde{T}$-action

$$t \cdot (f \otimes \alpha) := t \cdot f \otimes \chi^{\widetilde{u}}(t)\alpha,$$

where the $\widetilde{T}$-action on $\mathcal{O}_U$ is the usual one, given by $(t \cdot f)(x) = f(t^{-1} \cdot x)$. Let $\kappa\colon T \to \widetilde{T}$ be a homomorphism with $\pi \circ \kappa = \mathrm{id}_T$. Then $T$ acts via $\kappa$ on $\mathcal{O}_U \otimes \mathbb{C}(\widetilde{u})$ making

$$(p_*(\mathcal{O}_U \otimes \mathbb{C}(\widetilde{u})))^{H_0}$$

into a $T$-sheaf. We want to describe this $T$-sheaf in terms of certain divisors. To any $[\varrho, i] \in \Omega(\mathcal{S})$ with $\dim(\varrho) = 1$ we associate a Weil divisor on $X$ by setting

$$V_{[\varrho,i]} := \overline{T \cdot x_{[\varrho,i]}}.$$

**10.2 Lemma.** *For every $\widetilde{u} \in \widetilde{M}$, there is a canonical isomorphism of $T$-sheaves*

$$(p_*(\mathcal{O}_U \otimes \mathbb{C}(\widetilde{u})))^{H_0} \cong \mathcal{O}_X\left(\sum_{\substack{[\varrho,i] \in \Omega(\mathcal{S}),\\ \dim(\varrho)=1}} \langle \widetilde{u}, e_{[\varrho,i]}\rangle V_{[\varrho,i]}\right).$$

**Proof.** Let $L := \ker(P)$ and set

$$D := \sum_{\substack{[\varrho,i]\in\Omega(\mathcal{S}),\\ \dim(\varrho)=1}} \langle \widetilde{u}, e_{[\varrho,i]}\rangle V_{[\varrho,i]}, \qquad \widetilde{D} := \sum_{\substack{[\varrho,i]\in\Omega(\mathcal{S}),\\ \dim(\varrho)=1}} \langle \widetilde{u}, e_{[\varrho,i]}\rangle V_{\widetilde{\varrho}(i)}.$$

We claim that the natural map $\mathbb{C}(X) \to \mathbb{C}(U)$, $f \mapsto f \circ p$ defines an isomorphism of $T$-sheaves

$$\mathcal{O}_X(D) \to p_*(\mathcal{O}_U(\widetilde{D}))^{H_0}.$$

On the one hand, the sections of $\mathcal{O}_X(D)$ over a given affine chart $X_{[\sigma_{ii},i]} \subset X$ are generated by

$$\{\chi^u \in \mathrm{X}(T);\ \forall [\varrho,i] \prec [\sigma_{ii},i]\ \langle u, v_\varrho\rangle \geq -\langle \widetilde{u}, e_{[\varrho,i]}\rangle\}$$

On the other hand, the sections of $p_*(\mathcal{O}_U(\widetilde{D}))^{H_0}$ over $X_{[\sigma_{ii},i]}$ are generated by

$$\{\chi^{u'} \in \mathrm{X}(\widetilde{T});\ u'|_L = 0,\ \forall [\varrho,i] \prec [\sigma_{ii},i]\ \langle u', e_{[\varrho,i]}\rangle \geq -\langle \widetilde{u}, e_{[\varrho,i]}\rangle\}.$$

Since we have $\chi^u \circ p = \chi^{u \circ P}$ and both sheaves are coherent, our claim follows. Now the assertion follows from the fact that there is a further $T$-equivariant isomorphism, namely

$$p_*(\mathcal{O}_U \otimes \mathbb{C}(\widetilde{u}))^{H_0} \to p_*(\mathcal{O}_U(\widetilde{D}))^{H_0}, \qquad f \otimes \alpha \mapsto \alpha f \chi^{-\widetilde{u}}. \qquad \square$$

**Proof of Theorem 10.1.** First consider a coherent $H$-sheaf $\mathcal{F}$ on the toric prevariety $X = X_\mathcal{S}$, where $\mathcal{S}$ is an affine system of fans and $H$ acts on $X$ by means of a monomorphism $\varphi$ from $H$ to the acting torus $T$ of $X$.

Let $p\colon U \to X$ as described above, let $\widetilde{\varphi} := \kappa \circ \varphi \colon H \to \widetilde{T}$ and set $\widetilde{H} := \widetilde{\varphi}(H) \cdot H_0$. Since $\varphi$ is a monomorphism, $\mathcal{F}$ and hence also $p^*\mathcal{F}$ are in a canonical manner $\widetilde{H}$-sheaves. As in [5], p. 70, we obtain a finite-dimensional $\widetilde{H}$-module $V$ and exact sequences of $\widetilde{H}$-sheaves

$$\mathcal{O}_U \otimes V \longrightarrow p^*\mathcal{F} \longrightarrow 0, \qquad (p_*(\mathcal{O}_U \otimes V))^{H_0} \longrightarrow (p_*p^*\mathcal{F})^{H_0} \longrightarrow 0.$$

Using local resolutions of $\mathcal{F}$ by free sheaves, we see that the canonical morphism $\mathcal{F} \to (p_*p^*\mathcal{F})^{H_0}$ is in fact an isomorphism. Moreover, since the action of $\widetilde{H}$ on $V$ is diagonalizable, $\mathcal{O}_U \otimes V$ is a direct sum of $\widetilde{H}$-sheaves of the form $\mathcal{O}_U \otimes \mathbb{C}(\widetilde{u}_i)$, with some $\widetilde{u}_1, \ldots, \widetilde{u}_r \in \widetilde{M}$.

Now, consider the action of $H$ on the above $\widetilde{H}$-sheaves defined via $\widetilde{\varphi}\colon H \to \widetilde{H}$. Then Lemma 10.2 yields that the $H$-sheaf $(p_*(\mathcal{O}_U \otimes V))^{H_0}$ is equivariantly isomorphic to a direct sum of $H$-sheaves $\mathcal{O}_X(D_i)$. Hence we have an $H$-equivariant resolution of $\mathcal{F}$, given by

$$\bigoplus_{i=1}^r \mathcal{O}_X(D_i) \longrightarrow \mathcal{F} \longrightarrow 0.$$

In the general situation, notice that effectivity of the action of $H$ on $Y$ implies that $\varphi$ is injective. Now the assertion follows directly from applying the above considerations to the coherent $H$-sheaf $\imath_*\mathcal{F}$. $\square$

Jürgen Hausen, Fachbereich für Mathematik und Statistik, Universität Konstanz, D-78457 Konstanz, Germany. Email-adress: `Juergen.Hausen@uni-konstanz.de`